\def\b#1{{\bf #1}}   

\def\EPSpicture #1 by #2 (#3){
  \vbox to #2{
    \hrule width #1 height 0pt depth 0pt
    \vfill
    \special{illustration #3} 
    }
  }

\def\EPSscaledpicture #1 by #2 (#3 scaled #4){{
  \dimen0=#1 \dimen1=#2
  \divide\dimen0 by 1000 \multiply\dimen0 by #4
  \divide\dimen1 by 1000 \multiply\dimen1 by #4
  \EPSpicture \dimen0 by \dimen1 (#3 scaled #4)}
  }

\def\nref#1{(\ref{#1})}

\def\b#1{{\bf #1}}   

\def\horzpartitionline{\noalign{\vskip 5pt\hrule\vskip 5pt}}

\def\doubleuppervpartitionline#1{\smash{{\vrule height 8pt depth 6pt}\kern#1pt{\vrule height 8pt depth 6pt}}}             
\def\doublelowervpartitionline#1{\smash{{\vrule height 12pt depth 2pt}\kern#1pt{\vrule height 12pt depth 2pt}}}           
\def\doublemiddlevpartitionline#1{\smash{{\vrule height 12pt depth 6pt}\kern#1pt{\vrule height 12pt depth 6pt}}}          
\def\doublevariablevpartitionline#1#2#3{\smash{{\vrule height #1pt depth #2pt}\kern#3pt{\vrule height #1pt depth #2pt}}}  
\documentclass{article}[10]

\usepackage{graphicx}

\def\Width@semiexpanded{sx}

   {\hfill $\Box$                          
         \end{description}\end{flushleft}} 

\title{Algorithms, Initializations, and Convergence for the Nonnegative Matrix Factorization}

\author{Amy N. Langville$^\dagger$, Carl D. Meyer$^*$, Russell Albright$^\circ$, James Cox$^\circ$, and David Duling$^\circ$
}

\date{}

\begin{document}
\bibliographystyle{plain}

\maketitle

\begin{abstract}
It is well-known that good initializations can improve the speed and accuracy of the solutions of many nonnegative matrix factorization (NMF) algorithms \cite{Wild2003:thesis}.   Many NMF algorithms are sensitive with respect to the initialization of $\b W$ or $\b H$ or both.  This is especially true of algorithms of the alternating least squares (ALS) type \cite{Smildebook}, including the two new ALS algorithms that we present in this paper.  We compare the results of six initialization procedures (two standard and four new) on our ALS algorithms.  Lastly, we discuss the practical issue of choosing an appropriate convergence criterion.

\end{abstract}

\vskip .1in



\vskip .1in

\noindent
$^\dagger$
Department of Mathematics, \\
College of Charleston,\\
Charleston, SC 29424, USA \\
langvillea@cofc.edu \\
Phone: (843) 953-8021 \\ 
\newline
$^*$
Department of Mathematics and The Institute of Advanced Analytics, \\
North Carolina State University,
Raleigh, N.C. 27695-8205, USA  \\
meyer@math.ncsu.edu \\
Phone: (919) 515-2384 \\
Research supported in part by NSF CCR-ITR-0113121 and NSF DMS 9714811. \\
$^\circ$
SAS Institute, Inc., \\
Cary, NC 27513-2414, USA \\
\{russell.albright, james.cox, david.duling\}@sas.com \\
\newline

\newpage

\section{Introduction}
\label{Introduction}

Nonnegative data are pervasive.  Consider the following four important applications, each of which give rise to nonnegative data matrices.
\begin{itemize}
  \item In document collections, documents are stored as vectors.  Each element of a document vector is a count (possibly weighted) of the number of times a corresponding term appears in that document.  Stacking document vectors one after the other creates a nonnegative term-by-document matrix that represents the entire document collection numerically.
  \item Similarly, in image collections, each image is represented by a vector, and each element of the vector corresponds to a pixel.  The intensity and color of the pixel is given by a nonnegative number, thereby creating a nonnegative pixel-by-image matrix.
  \item For item sets or recommendation systems, the information for a purchase history of customers or ratings on a subset of items is stored in a non-negative sparse matrix.
  \item In gene expression analysis, gene-by-experiment matrices are formed from observing the gene sequences produced under various experimental conditions.
\end{itemize}
These are but four of the many interesting applications that create nonnegative data matrices (and tensors)
\cite{Langville2007:nmfsurvey}.

Three common goals in mining information from such matrices are: (1) to
automatically cluster similar items into groups, (2) to retrieve
items most similar to a user's query, and (3) identify interpretable critical dimensions within the  collection.  For the past decade, a technique called Latent
Semantic Indexing (LSI) \cite{BerryUSEbook}, originally conceived for the information retrieval problem  and later extended  to more general text mining problems, was a popular
means of achieving these goals.  LSI uses a well-known factorization
of the term-by-document matrix, thereby creating a low rank approximation of the original matrix.  This factorization, the singular value decomposition (SVD)
\cite{GolubVanLoanbook, MeyerSIAMbook}, is a classic technique in numerical
linear algebra.  

The SVD is easy to compute and works well for points (1) and (2) above, but not (3).  
The SVD does not provide users with any interpretation of its
mathematical factors or why it works so well. A common complaint from users is:
\emph{do the SVD factors reveal anything about the data collection?} 
Unfortunately, for the SVD, the answer to this question is no, as explained in the next section.  However, an
alternative and much newer matrix factorization, known as \emph{the
nonnegative matrix factorization (NMF), allows the question to be answered affirmatively.}   As a result, it can be shown that the NMF works nearly as well as the SVD on points (1) and (2), and further, can also achieve goal (3).

Most examples and applications of the NMF in this paper refer to text mining because this is the area with which we are most familiar.  However, the phrase ``term-by-document matrix"  which we will use frequently throughout this paper can just as easily be replaced with gene-by-observation matrix, purchase-by-user matrix, etc., depending on the application area.

\section{Low Rank Approximations}

Applications, such as text processing, data mining, and image
processing, store pertinent information in a huge matrix. 
This matrix $\b A$ is large, sparse, and often times nonnegative.  In
the last few decades, researchers realized that the data matrix
could be replaced with a related matrix, of much lower rank.  The
low rank approximation to the data matrix $\b A$ brought several
advantages.  The rank-$k$ approximation, denoted $\b A_k$, sometimes
required less storage than $\b A$.  But most importantly, the low
rank matrix seemed to give a much cleaner, more efficient
representation of the relationship between data elements.  The low
rank approximation identified the most essential components of the
data by ignoring inessential components attributed to noise,
pollution, or inconsistencies.  Several low rank
approximations are available for a given matrix: QR, URV, SVD, SDD, PCA, ICA, NMF, CUR, etc. \cite{OLeary1998, MeyerSIAMbook, Smildebook, Drineas2006}.  In this section, we
focus on two such approximations, the SVD and the NMF, that have been applied to data mining
problems.

\subsection{The Singular Value Decomposition}

In 1991, Susan Dumais \cite{Dumais1991} used the
singular value decomposition (SVD) to build a low rank
approximation to the term-by-document matrix of information retrieval.  In
fact, to build a rank-$k$ approximation
$\b A_k$ to the rank $r$ term-by-document matrix $\b A$, simply use the
$k$ most significant singular components, where $k < r$.  That is,
\begin{eqnarray*}
\b A_k = \sum_{i=1}^k \sigma_i \b u_i \b v_i^T = \b
U_k \b \Sigma_k \b V_k^T, 
\end{eqnarray*}
where $\sigma_i$ is the $i^{th}$ singular value of $\b A$, and $\b
u_i$ and $\b v_i^T$ are the corresponding singular vectors
\cite{GolubVanLoanbook}.  The technique of replacing $\b A$ with the truncated
$\b A_k$ is called Latent Semantic Indexing (LSI) because the low rank
approximation reveals meanings and connections between documents that
were hidden, or latent, in the original noisy data matrix $\b A$.  

Mathematically, the
truncated SVD has one particularly appealing property:  of all
possible rank-$k$ approximations, $\b A_k$ is the best approximation
in the sense that $\| \b A- \b A_k \|_F$ is as small as possible
\cite{BerryUSEbook, BerryJessup1999}.  Thus, the truncated SVD provides a nice baseline against which all other low-rank approximations can be judged for quantitative accuracy.
This optimality property is
also nice in practice.  Algorithms for computing the $k$ most
significant singular components are fast, accurate, well-defined, and robust 
\cite{Templatesbook, BerryUSEbook,
GolubVanLoanbook}.  Two different algorithms will produce the same results up to
roundoff error.  Such uniqueness and computational robustness are
comforting.  Another advantage of the truncated SVD concerns building successive
low rank approximations.  Once $\b A_{100}$ has been computed, 
no further computation is required if, for example, for sensitivity
analysis or comparison purposes, other \emph{lower} rank
approximations are needed.  That is, once 
$\b A_{100}$ is available, then $\b A_k$ is available for any $k \le
100$.

LSI and the truncated SVD dominated text mining research in the 1990s
\cite{Baeza-Yatesbook, BerryCIRbook, BerryUSEbook,
BerryFierro1996, BerryJessup1999, BerryOBrien1998,
Blomthesis, BlomRuhe2001, Kumar2005, Dumais1991, HugheyBerry2000,
Littman2000, JiangBerry2000, LetscheBerry1997, WittenBerry1998, Zha1998:kvalue,
Berry2001:levelsets}.  However, LSI is not perfect. 
For instance, while it first appeared that the low rank approximation
$\b A_k$ would save storage over the original matrix $\b A$,
experiments showed that this was not the case.  $\b A$ is generally
very sparse for text mining problems because only a small subset of
the terms in the collection are used in any particular document.  
No matter how sparse the original term-by-document
matrix is, the truncated SVD produces singular components that are almost
always completely dense.  In many cases, $\b A_k$ can require more
(sometimes much more) storage than $\b A$.  

Furthermore, $\b A$ is
always a nonnegative matrix, yet the singular components are mixed in
sign.  
The SVD's loss of the nonnegative structure of the
term-by-document matrix means that the factors of the truncated SVD provide no
interpretability.  
To understand this statement, consider
a particular document vector, say, column 1 of $\b A$.  The truncated
SVD represents document 1, $\b A_1$, as
\begin{eqnarray*}
\b A_1 \approx  \sigma_1 v_{11} \pmatrix{ \vdots \cr \b u_1 \cr \vdots} + 
       \sigma_2 v_{12}  \pmatrix{\vdots \cr \b u_2 \cr \vdots} +
\cdots +
        \sigma_k v_{1k} \pmatrix{\vdots \cr \b u_k \cr \vdots },  
\end{eqnarray*}
which reveals that document 1 is a linear combination of the singular
vectors
$\b u_i$, also called the basis vectors.  The scalar weight
$\sigma_i v_{1i}$ represents the contribution of basis
vector
$i$ in document 1.  Unfortunately, the mixed signs in $\b u_i$ and
$\b v_i$ preclude interpretation.

Clearly, the interpretability issues with the SVD's
basis vectors are caused by the mixed signs in the singular vectors.  Thus,
researchers proposed an alternative low rank approximation that maintained the
nonnegative structure of the original term-by-document matrix.  
As a result, the nonnegative matrix factorization (NMF) was created
\cite{LeeSeung1999, PaateroTapper1994}.  The NMF replaces the role played by
the singular value decomposition (SVD).  Rather than factoring $\b A$ as $\b U_k
\b \Sigma_k \b V_k^T$, the NMF factors $\b A$ as $\b W_k \b H_k$, where $\b W_k$
and $\b H_k$ are nonnegative.

\subsection{The Nonnegative Matrix Factorization}

Recently, the nonnegative matrix factorization (NMF) has been used to create a low rank approximation to $\b A$ that contains nonnegative factors called $\b W$ and $\b H$.
The NMF of a data matrix $\b A$ is created by solving the
following nonlinear optimization problem. 
\vskip -15pt
\begin{eqnarray}
\min \| \b A_{m \times n} &-& \b W_{m \times k} \b H_{k \times n} \|^2_F, \\
s.t.  &&  \b W \ge \b 0,\nonumber \\   && \b H \ge \b 0 \nonumber.
\label{NMFopt} 
\end{eqnarray}
The Frobenius norm is often used to measure the error between the original matrix $\b A$ and its low rank approximation $\b W \b H$, but there are other possibilities \cite{Dhillon2005, LeeSeung1999, Paatero1997}.
The rank of the approximation, $k$, is a parameter that must be set by the user.

The NMF is used in place of other low rank factorizations, such as the singular value decomposition (SVD) \cite{MeyerSIAMbook}, because of its two primary advantages: storage and interpretability. Due to the nonnegativity constraints, the NMF produces a so-called ``additive parts-based" representation \cite{LeeSeung1999} of the data. One consequence of this is that the factors $\b W$ and $\b H$ are generally naturally sparse, thereby saving a great deal of storage when compared with the SVD's dense factors.  

\begin{figure}[h!]
\begin{center}
\includegraphics[scale=.5]{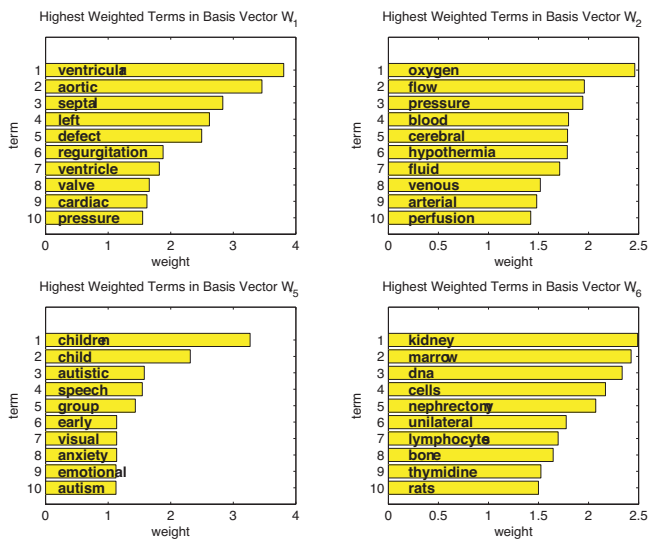}
\caption{Interpretation of NMF basis vectors on \texttt{medlars} dataset}
\label{figure:interpretdim}
\end{center}
\end{figure}

The NMF also has impressive benefits in terms of interpretation of its factors, which is, again, a consequence of the nonnegativity constraints. For example, consider a text processing application that requires the factorization of a term-by-document matrix $\b A_{m \times n}$.
In this case, $k$ can be considered the number of (hidden) topics 
present in the document collection.  In this case, $\b W_{m
\times k}$ becomes a term-by-topic matrix whose columns are the NMF
basis vectors.  The nonzero elements of column 1 of $\b W$ (denoted $\b W_1$),
which is sparse and nonnegative,
correspond to particular terms.  By considering the highest weighted terms in
this vector, one can assign a label or topic to basis vector 1.  Figure
\ref{figure:interpretdim} shows four basis vectors for one
particular term-by-document matrix, the \texttt{medlars} dataset of medical abstracts, available at \verb"http://www.cs.utk.edu/~lsi/".  
For those familiar with the domain of this dataset, the NMF allows users the ability to
interpret the basis vectors.  For instance, a user might attach
the label ``heart" to basis vector $\b W_1$ of Figure
\ref{figure:interpretdim}.
Similar interpretation holds for the
other factor $\b H$.  $\b H_{k \times n}$ becomes a topic-by-document
matrix with sparse nonnegative columns.  Element $j$
of column 1 of $\b H$ measures the strength to which topic $j$ appears
in document 1.

\begin{figure}[h!]
\begin{center}
\includegraphics[scale=1]{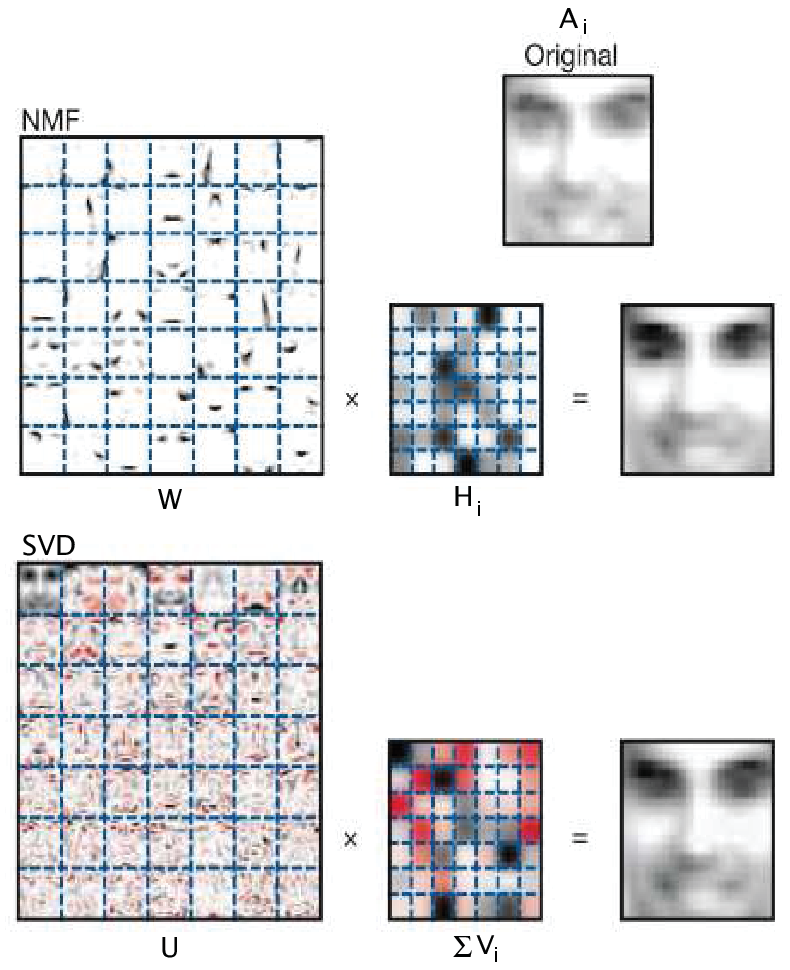}
\caption{Interpretation of NMF and SVD basis vectors on face dataset, from \cite{LeeSeung1999}}
\label{figure:LSnmffigures1}
\end{center}
\end{figure}

Another fascinating application of the NMF is image processing.  Figure \ref{figure:LSnmffigures1} clearly demonstrates two advantages of the NMF over the SVD. First, notice that the NMF basis vectors, represented as individual blocks in the $\b W$ matrix, are very sparse (i.e., there is much white space).  Similarly, the weights, represented as individual blocks in the $\b H_i$ vector, are also sparse.  On the other hand, the SVD factors are nearly completely dense.  Second, the basis vectors of the NMF, in the $\b W$ matrix, have a nice interpretation, as individual components of the structure of the face---ears, noses, mouths, hairlines.  The SVD basis vectors do not create an additive parts-based representation.
In addition, the gains in storage and interpretability do not come at a loss in performance. The NMF and the SVD perform equally well in reconstructing an approximation to the original image.

Of course, the NMF has its disadvantages too.  Other popular factorizations, especially the SVD, have strengths concerning uniqueness and robust
computation. Yet these become problems for the NMF.  There is no unique global
minimum for the NMF. The optimization problem of Equation
\nref{NMFopt} is convex in either $\b W$ or $\b H$, but not in
 both $\b W$ and
$\b H$, which means that the algorithms can only, if at all,
guarantee convergence to a local minimum.  In practice, NMF users
often compare the local minima from several different starting
points, using the results of the best local minimum found.  However,
this is prohibitive on large, realistically-sized problems.  Not
only will different NMF algorithms (and there are many now \cite{Langville2007:nmfsurvey})
produce different NMF factors, the
same NMF algorithm, run with slightly different parameters, can
produce different NMF factors.

\subsection{Summary of SVD vs. NMF}

We pause to briefly summarize the advantages of these two competing low rank approximations.
The properties and advantages of the SVD include: (1) an optimality property; the truncated SVD produces the best  rank-$k$ approximation (in terms of squared distances), (2) speedy and robust computation, (3) unique factorization; initialization does not affect SVD algorithms, and ( 4) orthogonality; resulting basis vectors are orthogonal and allow conceptualization of original data as vectors in space.   On the other hand, the advantages of NMF are: (1) sparsity and nonnegativity; the factorization maintains these properties of the original matrix,  (2) reduction in storage; the factors are sparse, which also results in easier application to new data, and (3) interpretability; the basis vectors naturally correspond to conceptual properties of the data.   

The strengths of one approximation become the weaknesses of another.  The most severe weakness of the NMF are its convergence issues.  Unlike the SVD and its unique factorization, there is no unique NMF factorization.  Because different NMF algorithms can converge to different local minima (and even this convergence to local minima is not guaranteed), initialization of the algorithm becomes critical.
In practice, knowledge of the application area can help guide initialization choices.  We will return to such initialization issues in Section \ref{Initializations}.

\section{ALS Algorithms for the NMF}
\label{ALSAlg}

Several currently popular NMF algorithms
\cite{ChuPlemmons2004, LeeSeung1999, LeeSeung2000, PaateroTapper1994,
Paatero1997, XuLiu2003} do not create sparse factors, which are desired for 
storage, accuracy, and interpretability reasons.  Even with adjustments to create sparse factors,
the improved algorithms
\cite{Hoyer2002, LeeSeung2000, PiperPlemmons2005, BerryPlemmons2004} exhibit an
undesirable
\emph{locking} phenomenon, as explained below.
Thus, in this section, we propose two new NMF algorithms \cite{SASSlides4}, called ACLS and AHCLS, that
produce sparse factors and avoid the so-called locking problem.

 Both algorithms are modifications to the
simple Alternating Least Squares (ALS) algorithm
\cite{PaateroTapper1994}, wherein $\b W$ is fixed and
$\b H$ is computed using least squares, then $\b H$ is fixed and
$\b W$ is computed using least squares, and so on, in alternating
fashion.  The method of alternating variables is a well-known
technique in optimization
\cite{Nocedal1999}.  
One problem with the first ALS algorithm applied to the NMF
problem (done by Paatero and Tapper in 1994 \cite{PaateroTapper1994}) was the lack of sparsity
restrictions.  To address this, the ACLS algorithm adds a reward for sparse
factors of the NMF.  The user sets the two parameters $\lambda_H$ and
$\lambda_W$ to positive values.  Increasing these values increases
the sparsity of the two NMF factors.  However, because there are no
upperbounds on these parameters, a user must resort to trial and error to find
the best values for $\lambda_H$ and
$\lambda_W$.  The more advanced AHCLS \cite{SASSlides4}, presented in Section \ref{AHCLS}, provides better
sparsity parameters with more intuitive bounds.

\subsection{The ACLS Algorithm}
\label{ACLS}

The ACLS (Alternating Constrained Least Squares) algorithm is implemented
differently than the original ALS algorithm \cite{PaateroTapper1994} because
issues arise at each alternating step, where a constrained least squares problem
of the following form 
\begin{eqnarray}
 \min_{\b h_{j}}
\|
\b a_{j} -
\b W \b h_{j}
\|_2^2 + \lambda_H \| \b h_{j} \|_2^2 \quad s.t. \quad \lambda_H \geq 0, \b
h_{j}
\geq
\b 0
\label{vectorACLS}
\end{eqnarray}
must be
solved.  The vectors $\b a_j$ and $\b h_j$ are columns of $\b A$ and $\b H$, respectively.  Notice that the decision variable $\b h_j$ must be nonnegative. There are
algorithms specifically designed for this nonnegative constrained least squares
problem. In fact, the NNLS algorithm of Lawson and Hanson
\cite{Bjorck1996, LawsonHansonLSbook} is so common that it appears as a
built-in function in MATLAB.  Unfortunately, the NNLS algorithm is
very slow, as it is an ``active set" method, meaning it can swap only one variable from the basis at a
time.  Even the faster version of the NNLS algorithm by Bro and de Jong \cite{Bro1997} is still not fast enough, and the NNLS step remains the computational bottleneck.
As a result, in practice, compromises are made. 
For example, a standard (unconstrained) least squares step is run \cite{Bjorck1996} and all
negative elements in the solution vector are set to 0.  This ad-hoc
enforcement of nonnegativity, while not theoretically appealing,
works quite well in practice.  The practical ACLS algorithm is shown
below.

\bigskip
\noindent
\framebox{
 \parbox[t]{6.3in}{
\centerline {\sc Practical ACLS Algorithm for
NMF} 

\smallskip
\noindent
input $\lambda_W$, $\lambda_H$ \\
$\b W$ = rand(m,k); \quad {\small \% initialize $\b W$ as
random dense matrix or use another initialization from Section
\ref{Initializations}}\\ 
for i = 1 : maxiter\\
\hspace*{.3in} {\sc (cls) } \quad \quad \ \ Solve for $\b H$ in matrix equation
$(\b W^T
\b W + \lambda_H \b I)\ \b H = \b W^T
\b A$. \quad   \ {\small \% for $\b W$ fixed, find $\b
H$}\\
\hspace*{.3in} {\sc (nonneg) } \, Set all negative elements in $\b H$ to 0.
\\
\hspace*{.3in} {\sc (cls) } \quad \quad \ \ Solve for $\b W$ in matrix equation
$(\b H
\b H^T  + \lambda_W \b I)\ \b W^T = \b H\b A^T
$. \ \ {\small \% for $\b H$ fixed, find $\b W$}\\
\hspace*{.3in} {\sc (nonneg) }  \,  Set all negative elements in $\b W$ to 0.
\\ end

}
}

\bigskip

\subsection{The AHCLS Algorithm}
\label{AHCLS}

ACLS uses a crude measure $\| \b x \|_2^2$ to approximate the sparsity of a vector $\b x$.
The AHCLS replaces this with a more sophisticated measure, $spar(\b x)$, which was invented by Hoyer \cite{Hoyer2004}. 
\begin{eqnarray*}
spar(\b x_{n \times 1})=\;
\matrix{\sqrt{n} -
\|\b x\|_1 /
\|
\b x\|_2\cr
\horzpartitionline
\noalign{\vskip -2pt} \sqrt{n}-1}
\end{eqnarray*}
In AHCLS (Alternating Hoyer-Constrained Least Squares), the user defines two scalars $\alpha_W$ and $\alpha_H$ in addition to $\lambda_H$ and $\lambda_W$ of ACLS.  For AHCLS, the two additional scalars $0 \leq \alpha_W, \alpha_H \leq 1$ represent a user's desired sparsity in each column of the factors.  These scalars, because they range from 0 to 1, match nicely with a user's notion of sparsity as a percentage.  Recall that $0 \leq \lambda_W, \lambda_H \leq \infty$ are positive weights associated with the penalties assigned to the density of $\b W$ and $\b H$.  Thus, in AHCLS, they measure 
how important it is to the user that $spar(\b W_{j*})=\alpha_W$ and
$spar(\b H_{j*})=\alpha_H$. Our experiments show that AHCLS does a better job of enforcing sparsity than ACLS does.  And the four AHCLS parameters are easier to set. For example, as a guideline, we recommend $0 \le
\lambda_W,\lambda_H
\le 1$, with of course, $0 \leq \alpha_W, \alpha_H \leq 1$.
The practical AHCLS algorithm, using matrix systems and ad-hoc enforcement of negativity, is below.  $\b E$ is the matrix of all ones.

\bigskip
\noindent
\framebox{
 \parbox[t]{6.3in}{
\centerline {\sc Practical AHCLS Algorithm for
NMF} 

\smallskip
\noindent
input $\lambda_W$, $\lambda_H$, $\alpha_W$, $\alpha_H$
$\b W$ = rand(m,k); \quad {\small \% initialize $\b W$ as
random dense matrix or use another initialization from Section
\ref{Initializations}}\\ 
$\beta_H = ((1-\alpha_H)\sqrt{k}+\alpha_H)^2$ \\
$\beta_W = ((1-\alpha_W)\sqrt{k}+\alpha_W)^2$ \\
for i = 1 : maxiter\\
\hspace*{.3in} {\sc (hcls) } \quad  \ \ \ Solve for $\b H$ in matrix equation
$(\b W^T
\b W + \lambda_H \beta_H \b I - \lambda_H \b E)\ \b H = \b W^T
\b A$. \\
\hspace*{.3in} {\sc (nonneg) } \, Set all negative elements in $\b H$ to 0.
\\
\hspace*{.3in} {\sc (hcls) } \quad  \ \ \ Solve for $\b W$ in matrix equation
$(\b H
\b H^T  + \lambda_W \beta_W \b I-\lambda_W \b E)\ \b W^T = \b H\b A^T
$.  \\
\hspace*{.3in} {\sc (nonneg) }  \,  Set all negative elements in $\b W$ to 0.
\\ end

}
}

\bigskip

\subsection{Advantages and Disadvantages of ACLS and AHCLS}
\label{AdDisadALS}

\subsubsection{Speed}

These algorithms have many advantages.  For
instance, rather than computing the vectors in $\b H$ column by column (as is done in \cite{BerryPlemmons2004}), thereby solving sequential least squares problems of the form of Equation (\ref{vectorACLS}), one matrix system solve can be executed.  Further, 
because each CLS step solves a \emph{small} $k \times k$ matrix system, ACLS
and AHCLS are the fastest NMF algorithms available (and faster than current
truncated SVD algorithms).  See Section \ref{ExpALS} for comparative run times.  They converge quickly and give
very accurate NMF factors.

\subsubsection{Sparsity}

Only $\b W$ must be initialized, and sparsity is
incorporated for both NMF factors.  We believe that avoidance of the
so-called
\emph{locking} phenomenon is one reason why the class of ALS algorithms works
well in practice.  Nearly all other NMF algorithms, especially those of the
multiplicative update class
\cite{Hoyer2002, Hoyer2004, LeeSeung1999, LeeSeung2000, BerryPlemmons2004:siam,
Shahnaz2004, BerryPlemmons2004}, lock elements when they become 0. That is, during the iterative
process, once an element in either $\b W$ or
$\b H$ becomes 0, it must remain 0.  For the basis vectors in the text mining problem, which are
stored in $\b W$, this means that in order to improve the objective
function, the algorithm can only remove terms from, not add terms to,
topic basis vectors.  As a result, once the algorithm starts
down a path toward a particular topic vector, it must continue in that
direction.  On the other hand, ALS algorithms do not {\it lock} elements, and
thus provide greater flexibility, allowing them to escape from a path
heading towards a poor local minimum.


\subsubsection{Convergence}

It has been proven that ALS algorithms will converge to a fixed point, but this fixed point may be a local extrema or a saddle point \cite{{Finesso2004, Gonzalez2005, Lin2005}}.  The ACLS and AHCLS algorithms with properly enforced nonnegativity, for example, by the NNLS algorithm, are known to converge to a local minimum \cite{Dhillon2005, Lin2005}. However, our ad-hoc enforcement of nonnegativity, which drastically speeds up the algorithm (and improves sparsity), means there are no proofs claiming convergence to a local minimum; saddle points are now possible.  (Actually, this is not so damning for our two ALS algorithms because most NMF algorithms suffer this same problem.  The few NMF algorithms believed to guarantee convergence to a local minimum have been proven otherwise \cite{{Finesso2004, Gonzalez2005, Lin2005}}.)
Our experiments 
\cite{SASSlides4, SASSlides3} and others \cite{Paatero1996, Paatero1999, PaateroTapper1994,
Paatero1997, Plemmons2005} have shown that the ALS fixed points can be superior to the results of other NMF algorithms.   

\subsubsection{Nonnegativity}

Clearly, ad-hoc enforcement of nonnegativity is theoretically unattractive. 
There are some alternatives to this ad-hoc enforcement of nonnegativity.  For instance, one could convert from an alternating least squares approach to an alternating linear programming approach, whereby nonnegativity of variables is enforced in a natural way by the simple constraints of the linear programming
formulation. Yet, this has the same problem as the NNLS algorithm, lengthy execution time.
A second alternative to ad-hoc enforcement of nonnegativity is to add negativity penalties in the form of logarithmic functions to the NMF objective function \cite{LuWu2004}.  This is a focus of future work.

\subsection{Numerical Experiments}
\label{ExpALS}

\begin{figure}[htbp]
\hspace*{-.2in}
\includegraphics[scale=1]{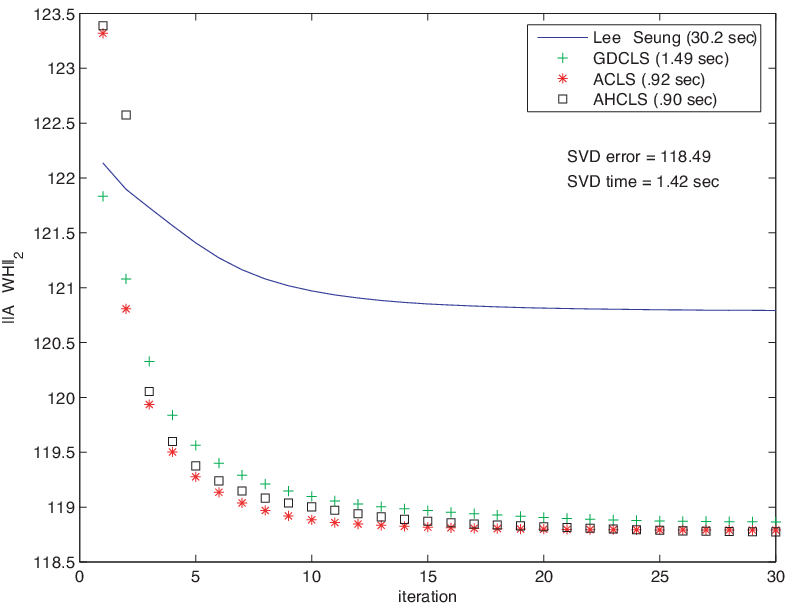}
\hfill
\includegraphics[scale=1]{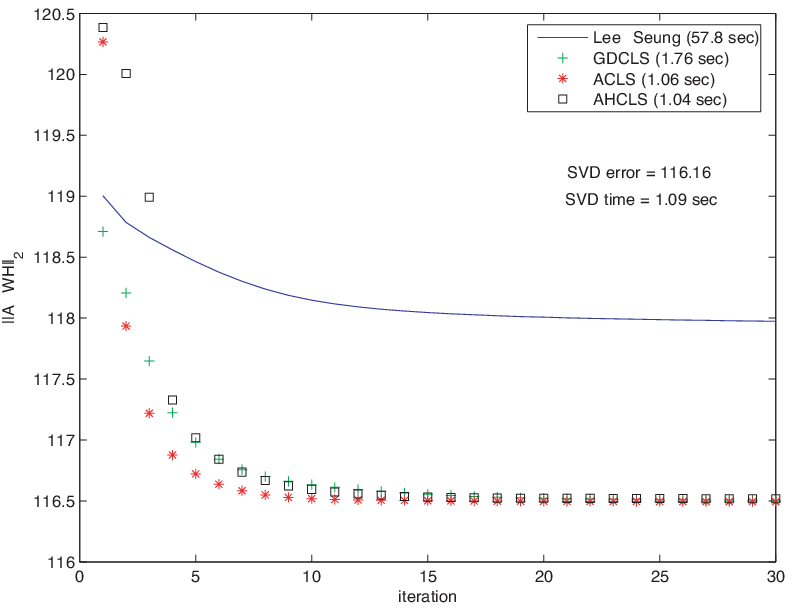}
\caption{Accuracy and Run-times of NMF Algorithms on \texttt{medlars} (left) and \texttt{cisi} (right) datasets}
\label{figure:NMFalgs}
\end{figure}

Figure \ref{figure:NMFalgs} compares our ACLS and AHCLS algorithms with the popular Lee-Seung mean squared error algorithm \cite{LeeSeung1999} and the GDCLS algorithm \cite{BerryPlemmons2004}.   We use our own implementation of GDCLS, which is much faster than the implementation presented in \cite{BerryPlemmons2004}.  The speed improvement results from our use of one matrix system rather than serial vector systems to solve the CLS step.  This implementation trick was described above for the ACLS and AHCLS algorithms.

To create Figure \ref{figure:NMFalgs}, we used the \texttt{medlars} dataset of medical abstracts and the \texttt{cisi} dataset of library science abstracts.
These figures clearly show how the ALS-type algorithms outperform the Lee-Seung multiplicative update algorithms in terms of accuracy and speed.  While the ALS-type algorithms provide similar accuracy to the GDCLS algorithm, they are much faster.  This speed advantage continues to hold for much larger collections like the \texttt{reuters10} collection, used in Section \ref{Initializations}.  (Details on the \texttt{reuters10} dataset appear in Section \ref{reuters10}.) On average the ACLS and AHCLS algorithms require roughly 2/3 the time of the GDCLS algorithm.  Figure \ref{figure:NMFalgs} also reports the error in the optimal rank-10 approximation required by the SVD.  Notice how close all NMF algorithms come to the optimal factorization error.  Also, notice that ACLS and AHCLS require less time than the SVD to produce such good, sparse, nonnegative factorizations.

 \section{Initializations}
\label{Initializations}

All NMF algorithms are iterative and it is well-known that they are sensitive to the initialization of $\b W$ and $\b H$ \cite{Wild2003:thesis}.
Some algorithms require that both $\b W$ and $\b H$ be
initialized
\cite{Hoyer2002, Hoyer2004, LeeSeung1999,
LeeSeung2000, PiperPlemmons2005}, while others require
initialization of only $\b W$
\cite{PaateroTapper1994, Paatero1997, Shahnaz2004, BerryPlemmons2004}.  
\emph{In all
cases, a good initialization can improve the speed and accuracy of the
algorithms, as it can produce faster convergence to an improved local minimum} \cite{Smildebook}. A good initialization 
can sidestep some of the convergence problems mentioned above, which is precisely why they are so important.
In this section, we compare several initialization procedures (two old and four new) by testing them on the ALS algorithms presented in Section \ref{ALSAlg}.   We choose to use the ACLS and AHCLS algorithms because they produce sparse accurate factors and require about the same time as the SVD.  Most other NMF algorithms require much more time than the SVD, often times orders of magnitude more time.

\subsection{Two Existing Initializations}
Nearly all NMF algorithms use 
simple \emph{random initialization}, i.e., $\b W$ and $\b H$ are initialized as
dense matrices of random numbers between 0 and 1.  It is well-known that random initialization does not generally provide a good first estimate for NMF algorithms \cite{Smildebook}, especially those of the ALS-type of \cite{Burdick1990, Li1993, Sanchez1990, Sands1980}.  Wild et
al.
\cite{Wild2003:thesis, Wild2003, Wild2004} have shown that the \emph{centroid initialization}, built from the centroid decomposition \cite{Dhillon2001} is a much better
alternative to random initialization.  Unfortunately, this decomposition is expensive as a preprocessing
step for the NMF.   Another advantage of ALS algorithms, such as our ACLS and AHCLS, is that they only require initialization of $\b W$.  In ALS algorithms, once $\b W^{(0)}$ is known, $\b H^{(0)}$ is computed quickly by a least squares computation.  As a result, we only discuss techniques for computing a good $\b W^{(0)}$.  

\subsection{Four New Initializations}
Some text mining software produces the SVD factors for other text tasks. Thus, in the event that the SVD factor $\b V$ is available, we propose a
\emph{SVD-centroid initialization} \cite{SASSlides3}, which initializes $\b
W$ with a centroid decomposition of the low dimensional SVD factor
$\b V_{n
\times k}$ \cite{Skillicorn2003}.  While the centroid decomposition of $\b A_{m \times n}$ can be too time-consuming, the centroid decomposition of $\b V$ is fast because $\b V_{n \times k}$ is much smaller than $\b A_{m \times n}$.
When the SVD factors are not available, we propose
a very inexpensive procedure called \emph{random Acol
initialization}.   Random Acol forms an initialization of each column
of the basis matrix $\b W$ by averaging
$p$ random columns of $\b A$.  It makes more sense to
build basis vectors from the given data, the sparse document vectors
themselves, than to form completely dense random basis vectors, as random
initialization does.  Random Acol initialization is very inexpensive, and lies
between random initialization and centroid initialization in terms of performance
\cite{SASSlides4, SASSlides3}.    

We also present two more initialization ideas, one inspired by the $\b C$ matrix of the CUR decomposition \cite{Drineas2006}, and another by the term co-occurrence matrix \cite{Sandler2005}.  We refer to these last two methods as \emph{random $\b C$ initialization} and \emph{co-occurrence initialization}, respectively.  The random $\b C$ initialization is similar to the random Acol method, except it chooses $p$ columns at random from the longest (in the 2-norm) columns of $\b A$, which generally means the densest columns since our text matrices are so sparse.  The idea is that these might be more likely to be the centroid centers.  The co-occurrence method first forms a term co-occurrence matrix $\b C = \b A \b A^T$. Next, the method for forming the columns of $\b W^{(0)}$ described as Algorithm 2 of \cite{Sandler2005} is applied to $\b C$.  The co-occurrence method is very expensive for two reasons. First, for text mining datasets, such as \texttt{reuters10}, $m>>n$, which means $\b C=\b A \b A^T$ is very large and often very dense too.  Second, the algorithm of \cite{Sandler2005} for finding $\b W^{(0)}$ is extremely expensive, making this method impractical.
All six initialization methods are summarized in Table \ref{table:InitMethods}. The two existing  methods appear first, followed by our four suggested methods.

\begin{table*}
\centering
\caption{Initialization Methods for the NMF}
\begin{center}
\begin{tabular}{|l|l|l|l|}
\hline
Name & Proposed by & Pros & Cons \cr
\hline
Random & Lee, Seung \cite{LeeSeung2000} & easy, cheap to compute & dense matrices, no intuitive basis \cr
Centroid & Wild et al. \cite{Wild2003:thesis} & reduces  \# NMF iterations,  & expensive, must run clustering \cr
   &  &  $\;\;\;\;$  firm, intuitive foundation  &  $\;\;\;\;$ algorithm on cols of $\b A$   \cr
   \hline
   SVD-Centroid & Langville \cite{SASSlides3} & inexpensive, reduces  \# NMF& SVD factor $\b V$ must be available \cr
   &  &  $\;\;\;\;$  iterations  & \cr
Random Acol & Langville \cite{SASSlides4} & cheap, sparse matrices built  & only slight decrease in number of   \cr
   &  &  $\;\;\;\;$ from original data & $\;\;\;\;$ NMF iterations \cr
 Random $\b C$ & Langville adapts & cheap, sparse & not very effective \cr
  & $\;\;\;\;$ from Drineas \cite{Drineas2006} &  & \cr
Co-occurrence & Langville adapts  & uses term-term similarities & large, dense co-occurrence matrix,\cr
   & $\;\;\;\;$ from Sandler \cite{Sandler2005}  &  $\;\;\;\;$     & $\;\;\;\;$  very expensive computation  \cr
  \hline
\end{tabular}
\end{center}
\label{table:InitMethods}
\end{table*}%

\subsection{Initialization Experiments with Reuters10 dataset} \label{reuters10}

The \texttt{reuters10} collection is our subset of the Reuters-21578 version of the Reuter's benchmark document collection of business newswire posts. The Reuters-21578 version contains over 20,000 documents categorized into 118 different categories, and is available online.\footnote{\texttt{\small http://www.daviddlewis.com/resources/testcollections/\\reuters21578/}} Our subset, the \texttt{reuters10} collection, is derived from the set of documents that have been classified into the top ten most frequently occurring categories. The collection contains 9248 documents from the training data of the ``ModApte split" (details of the split are also available at the website above).

The numbers reported in Table \ref{table:InitExps} were generated by applying the alternating constrained least squares (ACLS) algorithm of Section \ref{ALSAlg} with $\lambda_H=\lambda_W=.5$ to the \texttt{reuters10} dataset.
The error measure in this table is relative to the optimal rank-10 approximation given by the singular value decomposition.  For this dataset,  $\| \b A - \b U_{10} \b \Sigma_{10} \b V_{10}^T  \|_F=22656$.  Thus, for example, the error at iteration 10 is computed as
\begin{eqnarray*}
 \mbox{Error(10)}=  \frac{\| \b A - \b W^{(10)} \b H^{(10)} \|_F - 22656 }{22656}.
\end{eqnarray*}

\begin{table*}
\caption{Experiments with Initialization Methods for the NMF}
\begin{center}
\begin{tabular}{|l|c|c|c|c|c|c|c|}

\hline
Method & Time $\b W^{(0)}$  & Storage $\b W^{(0)}$ & Error(0) & Error(10) & Error(20) & Error(30) \cr
\hline
Random & .09 sec & 726K & 4.28\% & .28\% & .15\% & .15\% \cr
Centroid& 27.72 & 46K  & 2.02\% & .27\% & .18\% & .18\% \cr
   \hline
   SVD-Centroid & $.65^\dagger$ & 56K & 2.08\% & .06\% & .06\% & .06\% \cr
Random Acol$^*$ & .05 & 6K & 2.01\% & .21\% & .16\% & .15\% \cr
Random $\b C^{\circ}$ & .11 & 22K & 3.35\% & .29\% & .20\% & .19\% \cr
Co-occurrence & 3287  & 45K & 3.38\% & .37\% & .27\% & .25\% \cr
  \hline
  ACLS time & & & .37 sec & 3.42 & 6.78 & 10.29 \cr
  \hline
\end{tabular}
\end{center}
\label{table:InitExps}
$\dagger$ {\footnotesize provided $\b V$ of the SVD is already available} \\
\indent $^*$  {\footnotesize each column of $\b W^{(0)}$ formed by averaging 20 random columns of $\b A$} \\
\indent $^{\circ}$ {\footnotesize each column of $\b W^{(0)}$ formed by averaging 20 of the longest columns of $\b A$}\\
\end{table*}%

We distinguish between quantitative accuracy, as reported in Table \ref{table:InitExps}, and qualitative accuracy as reported in Tables \ref{table:BestBasisVectors} through \ref{table:CooccurrenceBasisVectors}.   For text mining applications, it is often not essential that the low rank approximation be terribly precise.  Often suboptimal solutions are ``good enough."  After reviewing Tables \ref{table:BestBasisVectors}--\ref{table:CooccurrenceBasisVectors}, it is easy to see why some initializations give better accuracy and converge more quickly.  They start with basis vectors in $\b W^{(0)}$ that are much closer to the best basis vectors found, as reported in Table \ref{table:BestBasisVectors}, which was generated by using the basis vectors associated with the best global minimum for the \texttt{reuters10} dataset, found by using 500 random restarts.  In fact, the relative error for this global minimum is .009\%, showing remarkable closeness to the optimal rank-10 approximation.  By comparing each subsequent table with Table \ref{table:BestBasisVectors}, it's clear why one initialization method is better than another.  The best method, SVD-centroid initialization, starts with basis vectors very close to the ``optimal" basis vectors of Table \ref{table:BestBasisVectors}.  On the other hand, random and random Acol initialization are truly random.   Nevertheless, random Acol does maintain one clear advantage over random initialization as it creates a very sparse $\b W^{(0)}$. The Random C and co-occurrence initializations suffer from lack of diversity.  Many of the longest documents in the \texttt{reuters10} collection appear to be on similar topics, thus, not allowing $\b W^{(0)}$ to cover many of the reuters topics.

Because the algorithms did not produce the ``wheat" vector always in column one of $\b W$, we have reordered the resulting basis vectors in order to make comparisons easier.  We also note that the nonnegative matrix factorization did produce basis vectors that cover 8 of the 10 ``correct" reuters classifications, which appear on the last line of Table \ref{table:BestBasisVectors}.  The two missing reuters classifications are \texttt{corn} and \texttt{grain}, both of which are lumped into the first basis vector labeled \texttt{wheat}.  This first basis vector does break into two separate vectors, one pertaining to \texttt{wheat and grain} and another to \texttt{corn} when the number of basis vectors is increased from $k=10$ to $k=12$.  We note that these categories have been notoriously difficult to classify, as previously reported in \cite{Dumais1998}.

\begin{table*}
\caption{Basis vectors of $\b W^{(50)}$ from \emph{Best Global Minimum} found for \texttt{reuters10}}
\begin{center}
\begin{tabular}{|c|c|c|c|c|c|c|c|c|c|}

\hline
$\b W_1^{(50)}$  & $\b W_2^{(50)}$ & $\b W_3^{(50)}$ & $\b W_4^{(50)}$ & $\b W_5^{(50)}$ & $\b W_6^{(50)}$   & $\b W_7^{(50)}$ & $\b W_8^{(50)}$ & $\b W_9^{(50)}$ & $\b W_{10}^{(50)}$  \cr
\hline

tonne & billion & share & stg & mln-mln & gulf & dollar & oil & loss & trade \cr
wheat & year & offer & bank & cts & iran & rate &  opec & profit & japan\cr
grain & earn & company & money & mln & attack & curr. & barrel & oper&  japanese\cr
crop & qrtr & stock & bill & shr & iranian  & bank &  bpd & exclude & tariff\cr
corn & rise & sharehol. & market & net & ship & yen & crude & net & import\cr
agricul.&  pct & common & england & avg & tanker & monetary & price & dlrs & reagan\cr
  \hline
  \texttt{wheat} & \texttt{earn} & \texttt{acquisition} &  & \texttt{interest} & \texttt{ship} & \texttt{frgn-exch.} & \texttt{oil} &   & \texttt{trade} \cr
 \hline
\end{tabular}
\end{center}
\label{table:BestBasisVectors}
\end{table*}%

\begin{table*}
\caption{Basis vectors of $\b W^{(0)}$  created by \emph{Random Initialization} for \texttt{reuters10}}
\begin{center}
\begin{tabular}{|c|c|c|c|c|c|c|c|c|c|}

\hline
$\b W_1^{(0)}$  & $\b W_2^{(0)}$ & $\b W_3^{(0)}$ & $\b W_4^{(0)}$ & $\b W_5^{(0)}$ & $\b W_6^{(0)}$   & $\b W_7^{(0)}$ & $\b W_8^{(0)}$ & $\b W_9^{(0)}$ & $\b W_{10}^{(0)}$  \cr
\hline
announce & wpp & formality & bulletin & matthews & dramatic & squibb & wag & cochran & erik \cr
medtec & reflagging & simply & awfully & nyt & boca raton & kuwaiti & oils & mln & support \cr
pac & kwik & moonie & blair & barrel & clever & dacca & hears & barriers & sale oil \cr
purina & tilbury & tmg & fresno & purina & billion & democrat & bwtr & deluxe & direct \cr
mezzanine & capacitor & bushnell & farm & june & bkne & induce & nestle & mkc & wheat \cr
foreign & grain & country & leutwiler & trend & clever & rate & fed. & econ. & aid \cr
  \hline
\end{tabular}
\end{center}
\label{table:RandomBasisVectors}
\end{table*}%

\begin{table*}
\caption{Basis vectors of $\b W^{(0)}$  created by \emph{Centroid Initialization} for \texttt{reuters10}}
\begin{center}
\begin{tabular}{|c|c|c|c|c|c|c|c|c|c|}

\hline
$\b W_1^{(0)}$  & $\b W_2^{(0)}$ & $\b W_3^{(0)}$ & $\b W_4^{(0)}$ & $\b W_5^{(0)}$ & $\b W_6^{(0)}$   & $\b W_7^{(0)}$ & $\b W_8^{(0)}$ & $\b W_9^{(0)}$ & $\b W_{10}^{(0)}$  \cr
\hline
tonne   &  bank & share & medar & cts & iran & rate & oil & stg & strike \cr
wheat   &   rate &    company & mdxr & mmln & gulf & dollar & trade & bill & port \cr
 grain     &  dollar &  offer & mlx  & loss & attack & bank & price & take-up & union \cr
corn      &   billion &  pct & mlxx  & net & iranian & currency & barrel & drain & seaman \cr
crop    &  pct & stock & mich & shr & missile & market & japan & mature & worker \cr
agriculture &  trade & dlrs & troy & dlrs & tanker & monetary & opec & money & ship \cr
  \hline
\end{tabular}
\end{center}
\label{table:CentroidBasisVectors}
\end{table*}%

\begin{table*}
\caption{Basis vectors of $\b W^{(0)}$  created by \emph{SVD-Centroid Initialization} for \texttt{reuters10}}
\begin{center}
\begin{tabular}{|c|c|c|c|c|c|c|c|c|c|}

\hline
$\b W_1^{(0)}$  & $\b W_2^{(0)}$ & $\b W_3^{(0)}$ & $\b W_4^{(0)}$ & $\b W_5^{(0)}$ & $\b W_6^{(0)}$   & $\b W_7^{(0)}$ & $\b W_8^{(0)}$ & $\b W_9^{(0)}$ & $\b W_{10}^{(0)}$  \cr
\hline
tonne & billion & share & bank & cts & iran & dollar & oil & loss & trade \cr
wheat & year & offer & money & shr & gulf & rate &  barrel & oper & japan\cr
grain & earn & company & rate & mln & attack & curr. & opec & profit &  japanese\cr
corn & qrtr & stock & stg & net & iranian  & yen &  crude & cts & tariff\cr
crop & rise & pct & market & mln-mln & missile & japan & bpd & mln & import\cr
agricul.&  pct & common & pct & rev & ship & economic & price & net & country \cr
  \hline
\end{tabular}
\end{center}
\label{table:SVDBasisVectors}
\end{table*}%

\begin{table*}
\caption{Basis vectors of $\b W^{(0)}$  created by \emph{Random Acol Initialization} for \texttt{reuters10}}
\begin{center}
\begin{tabular}{|c|c|c|c|c|c|c|c|c|c|}

\hline
$\b W_1^{(0)}$  & $\b W_2^{(0)}$ & $\b W_3^{(0)}$ & $\b W_4^{(0)}$ & $\b W_5^{(0)}$ & $\b W_6^{(0)}$   & $\b W_7^{(0)}$ & $\b W_8^{(0)}$ & $\b W_9^{(0)}$ & $\b W_{10}^{(0)}$  \cr
\hline
 mln & fee & agl & mln & mark &loss & official & dlrs & bank & trade \cr
 denman & mortg. & tmoc & dlrs &  mannes. & mln & piedmont &  oper & bancaire & viermetz \cr
 dlrs & billion & bank &share & dividend  & cts & dollar &billion & austral & mln \cr
 ecuador & winley & pct & seipp & mln  &maki  & interest & loss & neworld & nwa \cr
 venezuela & mln & company & billion      & dieter &name  & tokyo & texaco & datron & cts \cr
 revenue & fed & maki & dome & gpu & kato & japanese & pennzoil & share & builder \cr
  \hline
\end{tabular}
\end{center}
\label{table:RandomAcolBasisVectors}
\end{table*}%

\begin{table*}
\caption{Basis vectors of $\b W^{(0)}$  created by \emph{Random $\b C$ Initialization} for \texttt{reuters10}}
\begin{center}
\begin{tabular}{|c|c|c|c|c|c|c|c|c|c|}

\hline
$\b W_1^{(0)}$  & $\b W_2^{(0)}$ & $\b W_3^{(0)}$ & $\b W_4^{(0)}$ & $\b W_5^{(0)}$ & $\b W_6^{(0)}$   & $\b W_7^{(0)}$ & $\b W_8^{(0)}$ & $\b W_9^{(0)}$ & $\b W_{10}^{(0)}$  \cr
\hline
analyst & dollar & econ. & bank & market & analyst & analyst & analyst & trade & rate \cr
lawson & rate & policy & rate & bank &   market & industry & bank & dollar & trade \cr
market & econ. & pct & market & analyst & trade & price & currency & japan & official \cr
trade & mark & cost & currency & price &  pct & market  & japan & price & bank \cr
sterling & bank & growth & dollar & mark &   last & believe & billion & japanese & market \cr
dollar & rise & trade & trade & good &   official & last & cut & pct & econ. \cr
  \hline
\end{tabular}
\end{center}
\label{table:RandomCBasisVectors}
\end{table*}%

\begin{table*}
\caption{Basis vectors of $\b W^{(0)}$  created by \emph{Co-occurrence Initialization} for \texttt{reuters10}}
\begin{center}
\begin{tabular}{|c|c|c|c|c|c|c|c|c|c|}

\hline
$\b W_1^{(0)}$  & $\b W_2^{(0)}$ & $\b W_3^{(0)}$ & $\b W_4^{(0)}$ & $\b W_5^{(0)}$ & $\b W_6^{(0)}$   & $\b W_7^{(0)}$ & $\b W_8^{(0)}$ & $\b W_9^{(0)}$ & $\b W_{10}^{(0)}$  \cr
\hline
dept. & average & agricul. & national & farmer & rate-x & aver price & plywood & wash. & trade \cr
wheat & pct & wheat & bank & rate-x & natl & average &  aqtn & trade & japan\cr
agricul. & rate & tonne & rate & natl & avge & price & aequitron & japan &  billion\cr
tonne & price & grain & pct & avge & farmer  & yield &  medical & official & market \cr
usda & billion & farm & oil & cwt & cwt & billion & enzon & reagan & japanese \cr
corn &  oil & dept. & gov. & wheat & wheat & bill & enzon & pct & import \cr
  \hline
\end{tabular}
\end{center}
\label{table:CooccurrenceBasisVectors}
\end{table*}%

 \section{Convergence Criterion}
\label{ConvCriterion}

Nearly all NMF algorithms use the simplest possible convergence criterion, i.e., run for a fixed number of iterations, denoted \texttt{maxiter}.  This criterion is used so often because the natural criterion, stop when $\| \b A - \b W \b H \| \le \epsilon$, requires more expense than most users are willing to expend, even occasionally.
Notice that \texttt{maxiter} was the convergence criterion used in the ACLS and AHCLS algorithms of Section \ref{ALSAlg}. However, a fixed number of iterations is not a mathematically appealing way to control the number of iterations executed because the most appropriate value for \texttt{maxiter} is problem-dependent.

In this section, we use the ACLS algorithm applied to the \texttt{cisi} dataset to compare two convergence criterion: the natural but more expensive Frobenius norm measure, and our proposed angular measure, which we describe in the next paragraph.   We note that we used the most efficient implementation of the Frobenius measure \cite{Langville2007:nmfsurvey}, which exploits the trace form of the Frobenius norm.  
$$
\| \b A - \b W \b H \|_F^2 = trace(\b A^T \b A) - 2 \, trace(\b H^T \b W^T \b A) + trace(\b H^T \b W^T \b W \b H).
$$
In this equation $trace(\b A^T \b A)$ is a constant that does not change throughout the iterations, and thus, is only computed once and stored.  At each iteration $2 \, trace(\b H^T \b W^T \b A)$ and $trace(\b H^T \b W^T \b W \b H)$ must be computed.  However, some calculations required by these traces, such as $\b W^T \b A$ and $\b W^T \b W$, are already available from the least squares steps, and hence, need not be recomputed.

\begin{figure}[htbp]
\hspace*{-.2in}
\includegraphics[scale=1]{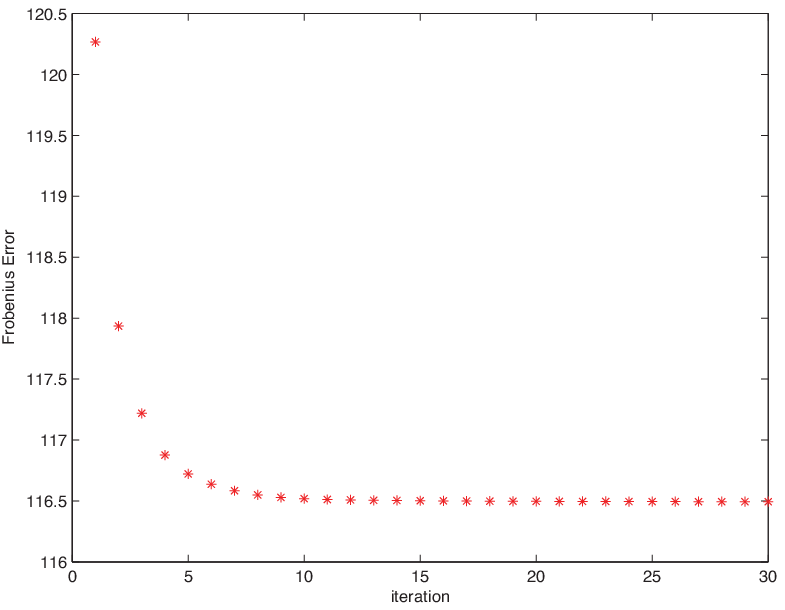}
\hfill
\includegraphics[scale=1]{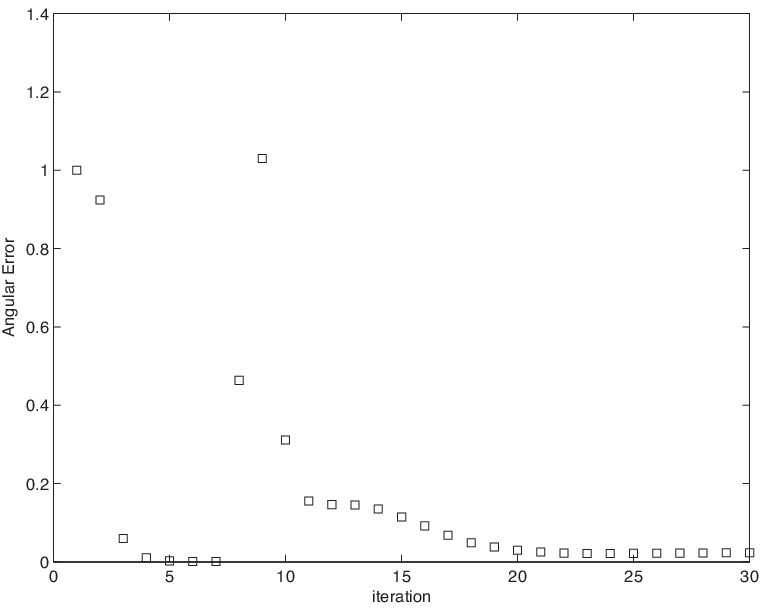}
\caption{Frobenius convergence measure (left) and angular convergence measure (right) of ACLS Algorithm on \texttt{cisi} datasets}
\label{figure:convcrit}
\end{figure}

Our angular convergence measure is moderately inexpensive in storage and computation, and is intuitively appealing.  Simply measure the angle $\theta_i$ between successive topic vectors, i.e., $\b W_i^{(j+1)}$ and $\b W_i^{(j)}$ at iterations $j$ and $j+1$. Once $\theta_i \le \epsilon$ for $i=1, \ldots, k$, stop because the topic vectors have converged satisfactorily.   Mitchell and Burdick \cite{Mitchell1994} have shown that, in a different context, a similar measure converges simultaneously with the expensive convergence criterion based on the objective function, $\| \b A - \b W \b H \|$ \cite{Smildebook}.  However, Figure \ref{figure:convcrit} clearly shows one problem with the angular convergence measure---it does not maintain continual descent, because the basis vectors compared from one iteration to the next are not required to maintain any fixed column order.  After several iterations, the column ordering of the basis vectors in $\b W$ is less likely to change, making the angular measure more useful in later iterations of the algorithm.   The angular convergence measure is much less expensive to compute than the Frobenius measure, but does require additional storage of the $\b W$ matrix from the previous iteration.
We note in practice, that regardless of the chosen convergence criterion, it is wise to only compute the measure every five or so iterations after some burn-in period.  

The recent 2005 reference by Lin \cite{Lin2005} mentioned the related convergence criterion issue of stationarity.  The fact that $\| \b A - \b W \b H \|$ (or some similar objective function) levels off does not guarantee stationarity.   Lin advocates a stationarity check once an algorithm has stopped.  For instance, the stationarity checks of Chu and Plemmons \cite{ChuPlemmons2004} may be used.  Lin \cite{Lin2005} proposes a convergence criterion, that simultaneously checks for stationarity, and fits nicely into his projected gradients algorithm.  We agree that a stationarity check should be conducted on termination.

\section{Conclusion}

The two new NMF algorithms presented in this paper, ACLS and AHCLS, are some of the fastest available, even faster than truncated SVD algorithms.  However, while the algorithms will converge to a stationary point, they cannot guarantee that this stationary point is a local minimum.  If a local minimum must be achieved, then we recommend using the results from a fast ALS-type algorithm as the initialization for one of the slow algorithms \cite{Lin2005} that do guarantee convergence to a local minimum.  In this paper we also presented several alternatives to the common, but poor, initialization technique of random initialization.  Lastly, we proposed an alternative stopping criterion that practical implementations of NMF code should consider.  The common stopping criterion of running for a fixed number of iterations should be replaced with a criterion that fits the context of the users and their data.  For many applications, iterating until $\| \b A - \b W \b H \|$ reaches some small level is unnecessary, especially in cases where one is most interested in the qualitative results produced by the vectors in $\b W$.  In such cases, our proposed angular convergence measure is more appropriate.

\bibliography{researchcite}

\end{document}